\documentclass{amsart}

\pdfoutput=1

\usepackage{amsmath,amsfonts,amsthm}
\usepackage{mathtools} % for abs and norm
\usepackage{MnSymbol} % for \llangle and \rrangle
\usepackage{enumerate}
\usepackage{tikz-cd}

\newcommand{\eps}{\varepsilon}
\newcommand{\Weyl}{\mathcal{W}}
\newcommand{\FFun}{\mathcal{F}}
\newcommand{\GFun}{\mathcal{G}}
\newcommand{\AFun}{A}
\newcommand{\catFields}{\textbf{\textup{Fields}}_{/k_0}}

\newcommand{\catRings}{\textbf{\textup{Rings}}}
\newcommand{\catSets}{\textbf{\textup{Sets}}}
\newcommand{\catAb}{\textbf{\textup{Ab}}}
\newcommand{\Ztwo}{\mathbb{Z}/2\mathbb{Z}}
\newcommand{\sw}{sw}
\newcommand{\swmod}{\widetilde{sw}}

\DeclareMathOperator{\sep}{sep}

\DeclareMathOperator{\Inv}{Inv}
\DeclareMathOperator{\NInv}{NInv}
\DeclareMathOperator{\Gal}{Gal}
\DeclareMathOperator{\Tr}{Tr}
\DeclareMathOperator{\Hom}{Hom}
\DeclareMathOperator{\res}{res}

\DeclarePairedDelimiter\floor{\lfloor}{\rfloor}

\newtheorem{theorem}{Theorem}[section]
\newtheorem*{theorem*}{Theorem}
\newtheorem{proposition}[theorem]{Proposition}
\newtheorem*{proposition*}{Proposition}

\newtheorem*{conjecture*}{Conjecture}

\newtheorem*{question*}{Question}
\newtheorem{lemma}[theorem]{Lemma}
\newtheorem*{lemma*}{Lemma}

\newtheorem*{corollary*}{Corollary}

\theoremstyle{definition}
\newtheorem{definition}[theorem]{Definition}
\newtheorem*{definition*}{Definition}

\newtheorem*{example*}{Example}
\newtheorem{remark}[theorem]{Remark}
\newtheorem*{remark*}{Remark}

\newtheorem*{fact*}{Fact}
\newtheorem*{note*}{Note}

\usepackage{hyperref}

\title{Witt invariants of Weyl groups}
\author{Tamar Lichter Blanks}
\date{}
\address{Department of Mathematics, Fordham University,
New York, NY 10023}
\email{tamarlichterblanks@gmail.com}
\thanks{\copyright \ 2024. This manuscript version is made available under the CC-BY-NC-ND 4.0 license \url{https://creativecommons.org/licenses/by-nc-nd/4.0/}.}

\begin{document}

\begin{abstract}
We describe the Witt invariants of a Weyl group over a field $k_0$ by giving generators for the $W(k_0)$-module of Witt invariants, under the assumption that the characteristic of $k_0$ does not divide the order of the group. For the Weyl groups of types $B_n$, $C_n$, $D_n$, and $G_2$, we show that the Witt invariants are generated as a $W(k_0)$-algebra by trace forms and their exterior powers, extending a result due to Serre in type $A_n$. Many of our computational methods are applicable to computing Witt invariants of any smooth linear algebraic group over $k_0$, including a technique for lifting module generators from cohomological invariants to Witt invariants.
\end{abstract}

\maketitle

\section{Introduction}

A Witt invariant is an invariant of algebraic objects over fields, valued in the Witt ring of quadratic forms. The trace form of a field extension, or more generally an {\' e}tale algebra, is a key example. In this paper we describe the Witt invariants that assign a class of quadratic forms to each $G$-Galois extension, for $G$ a Weyl group. 

Witt invariants were defined by J.-P. Serre \cite[\S 27.3, p\ 64]{Serre03_GMS} as part of the broader theory of invariants, as we discuss below and in Section \ref{subsec:introduction_invariants}. The precise definition of a Witt invariant is as follows. Fix a base field $k_0$ of characteristic not equal to $2$, and let $\catFields$ be the category of field extensions of $k_0$. For $\GFun$ a functor $\catFields \to \catSets$, a \emph{Witt invariant} of $\GFun$ is a natural transformation
\[
	w : \GFun(-) \to W(-),
\]
where $W(-)$ is the Witt ring functor; see also Definitions \ref{def:Invariant} and \ref{def:WittInvariant}. The intuition is that for each field extension $k$ of $k_0$, one has a set $\GFun(k)$ consisting of the algebraic objects over $k$ that one is interested in, and for each further extension $k_0 \subseteq k \subseteq k'$ one has an ``extension of scalars'' map $\GFun(k) \to \GFun(k')$. The invariant $w$ assigns an element of $W(k)$ to each element of $\GFun(k)$, and this assignment is compatible with extension of scalars. In many cases of interest, $\GFun$ can be expressed as a functor of the form $H^1(-, G)$, where $G$ is a smooth linear algebraic group over $k_0$ (see \cite[\S 2-3]{Serre03_GMS} for a nice list of examples). The Witt invariants of $H^1(-, G)$ are referred to as simply the Witt invariants of $G$. 

\medskip

The Weyl groups are a class of finite groups that appear in representation theory and are associated with Lie algebras (see for example \cite[Chapter III]{Humphreys1972}). In this paper we describe the Witt invariants of a Weyl group $G$, over a field $k_0$ of characteristic not dividing $|G|$, by giving generators for the $W(k_0)$-module of Witt invariants:~this is Theorem \ref{thm:WeylLiftGen}. For the Weyl groups of types $B_n$, $D_n$ and $G_2$, we show in Theorems \ref{thm:InvWeylBn}, \ref{thm:InvWeylDn}, and \ref{thm:InvWeylG2} respectively that the Witt invariants are generated as a $W(k_0)$-algebra by trace forms and their exterior powers. (The Weyl groups of types $B_n$ and $C_n$ are isomorphic, so we treat these cases simultaneously.) Many of our computational techniques, which are described in Section \ref{sec:Computations}, are applicable to computing Witt invariants of $G$ for any smooth linear algebraic group $G$ over $k_0$. 

Witt invariants belong to the more general theory of invariants defined by Serre \cite{GaribaldiMerkurjevSerre03}, as a related but potentially richer alternative to \emph{cohomological invariants} with coefficients in $\Ztwo$. As we discuss in Section \ref{subsec:introduction_invariants}, by the proof of the Milnor conjectures, one has a map from the degree $n$ Witt invariants of $G$ to the degree $n$ cohomological invariants of $G$. 

A challenge arises when one tries to translate information about the cohomological invariants of $G$ to information about the Witt invariants of $G$. The source of the problem is that there is no canonical way to lift an element of the graded ring associated to an ideal to the original filtered ring. To address this difficulty and obtain our results for Weyl groups, we establish sufficient conditions for lifting module generators from cohomological invariants to Witt invariants (Proposition \ref{prop:InvGenMod}), and combine this technique with an analysis of the relationships between particular invariants of note, including Stiefel-Whitney classes, exterior powers, and trace forms.

\subsubsection*{Acknowledgements} This article is an updated version of my PhD thesis, and I am deeply grateful to Daniel Krashen for his advice and mentorship. I thank Saurabh Gosavi for many valuable conversations, particularly those pertaining to Lemma \ref{lem:ProductInvar}; Stefan Gille for explaining the proof of Proposition \ref{prop:EquivGradedH}; Eva Bayer-Fluckiger, Lisa Carbone, Skip Garibaldi, Stefan Gille, David Harbater, Patrick McFaddin, Jean-Pierre Serre, and Charles Weibel for their insightful comments; and the anonymous referee for clear and excellent suggestions, in particular Lemma \ref{lem:Weyl_Char_Implies_BddInvar}. This material is based upon work supported by the National Science Foundation Graduate Research Fellowship under Grant No.~1842213.

\subsubsection*{Notation and conventions} All fields considered will be assumed to have characteristic\footnote{There are several reasons for the assumption on the characteristic, the most immediate being that we make use of the ring structure of $W(k)$ (not just the group structure). We also need this assumption for our application of the Milnor conjectures, and the relevant results on cohomological invariants in \cite{hirsch2020decomposability} are only available in characteristic not equal to $2$.} not equal to $2$. For a field $k$, let $k^{\sep}$ denote its separable closure and $\Gamma_k$ denote its absolute Galois group $\Gal(k^{\sep}/k)$. The notation $H^n(k, C)$ refers to degree $n$ Galois cohomology with coefficients in $C$. The categories of sets, abelian groups, and rings are denoted by $\catSets$, $\catAb$, and $\catRings$ respectively. The category $\catFields$ is the category of field extensions of $k_0$: the objects are fields $k$ with $k_0 \subseteq k$ and the morphisms are morphisms of fields $k \to k'$ that commute with the inclusions of $k_0$ into $k$ and $k'$. The Weyl groups of types $B_n$, $C_n$, $D_n$, and $G_2$ are denoted $\Weyl(B_n)$, $\Weyl(C_n)$, $\Weyl(D_n)$, and $\Weyl(G_2)$ respectively.

\section{Background}
\label{sec:intro_invariants_and_trace_forms}

In Section \ref{subsec:introduction_invariants}, we establish some basic facts about invariants, including the connection between Witt invariants and cohomological invariants via the Milnor conjectures \cite{Voevodsky2003, Morel2005Milnors, OrlovVishikVoevodsky2007}. 
Section \ref{subsec:traceform} is about trace form invariants in particular. We recall the result, due to Serre, that the Witt invariants of $S_n$ are generated by the trace form and its exterior powers.

\subsection{Invariants}
\label{subsec:introduction_invariants}

Fix a base field $k_0$. Following \cite{Serre03_GMS}, we let $G$ be a smooth linear algebraic group over $k_0$ (see \cite{springer1998linear}).

\begin{definition}\label{def:Invariant} 
Let $\AFun : \catFields \to \catAb$ be a functor. An \emph{$\AFun$-invariant of $G$} is a natural transformation
\[
	H^1(-, G) \to \AFun(-)
\]
of functors $\catFields \to \catSets$. The set of $\AFun$-invariants of $G$ is denoted $\Inv_{k_0}(G, \AFun)$.
\end{definition}

Note than an invariant $a \in \Inv_{k_0}(G, \AFun)$ is not required to send the trivial cocycle in $H^1(k, G)$ to $0 \in \AFun(k)$; in other words, $a$ need not be a natural transformation of functors from fields to pointed sets. If $a$ does send the trivial cocycle to $0$ for all fields $k/k_0$, then $a$ is said to be \emph{normalized}.

\begin{remark}\label{rmk:StructureOfInv}
The set $\Inv_{k_0}(G, \AFun)$ has the structure of an abelian group: for example, given invariants $a$ and $a'$, one defines $a+a'$ by
\[
	(a+a')_k : T \mapsto a_k(T) + a'_k(T)
\]
for every field extension $k$ of $k_0$ and every $T \in H^1(k, G)$. If $\AFun$ is valued in $\catRings$ then $\Inv_{k_0}(G, \AFun)$ is a ring as well. There is an injective group (or ring) homomorphism $\AFun(k_0) \hookrightarrow \Inv_{k_0}(G, \AFun)$, sending $x \in \AFun(k_0)$ to the \emph{constant invariant} at $x$, that is, the invariant sending every $T \in H^1(k, G)$ to the image $x_k$ of $x$ under $\AFun(k_0) \to \AFun(k)$. 
\end{remark}

In this paper we will be concerned with the functor $W(-)$ that sends each field $k$ to the Witt ring $W(k)$ and each field extension $k \subseteq k'$ to the extension of scalars map $W(k) \to W(k')$.

\begin{definition}\label{def:WittInvariant}
A \emph{Witt invariant} of $G$ is an element of $\Inv_{k_0}(G, W)$.
\end{definition}

By Remark \ref{rmk:StructureOfInv}, $\Inv_{k_0}(G, W)$ is a $W(k_0)$-algebra, so one can describe the Witt invariants of $G$ by finding generators for this algebra. Proposition \ref{prop:InvGenMod}, Theorem \ref{thm:WeylLiftGen}, and Theorems \ref{thm:InvWeylG2}, \ref{thm:InvWeylBn}, and \ref{thm:InvWeylDn} are all results of this type.

\medskip

While the definition of an $\AFun$-invariant is fairly general, previous work has emphasized the classification of \emph{cohomological invariants}. By definition, a cohomological invariant is an element of $\Inv_{k_0}(G, H)$, where $H(-) = \bigoplus_{n \geq 0} H^n(-, C)$ for some discrete $\Gamma_k$-module $C$. In \cite{GaribaldiMerkurjevSerre03}, S. Garibaldi, A. Merkurjev, and J.-P. Serre described the cohomological invariants of some particular classes of groups. For example, taking $C = \Ztwo$, Serre classified the cohomological invariants of $G = S_n$ and of $G = O_n$, corresponding by Galois descent to the cohomological invariants of {\' e}tale algebras and quadratic forms, respectively. More recently C. Hirsch \cite{hirsch2020decomposability} and Serre \cite{serre2018cohomological} described the cohomological invariants of Weyl groups for $C = \Ztwo$. (In fact the results of \cite{hirsch2020decomposability} are more general:~Hirsch classified the invariants $\Inv_{k_0}(G, M_*)$ where $G$ is a Weyl group, $M_*$ is a cycle module in the sense of \cite{Rost1996}, and $k_0$ is a field whose characteristic does not divide $|G|$. Note, however, that the Witt ring functor $W(-)$ is not a cycle module.) In this paper we will always take $H(-)$ to be the functor $\bigoplus_{n \geq 0} H^n(-, \Ztwo)$ and $H^n(-)$ to be the functor $H^n(-, \Ztwo)$. 

The connection between the Witt invariants of $G$ and the cohomological invariants of $G$ comes from the proof of the Milnor conjectures. Let $I^n(k)$ be the $n$th power of the fundamental ideal in $W(k)$. By the Milnor conjectures \cite{Voevodsky2003, Morel2005Milnors, OrlovVishikVoevodsky2007}, for each $k/k_0$ and $n \geq 0$ there is a well-defined, surjective group homomorphism
\[
	e_n(k) : I^n(k) \to H^n(k, \Ztwo)
\]
with kernel $I^{n+1}(k)$, sending the Pfister form $\llangle \alpha_1, \ldots, \alpha_n \rrangle = \langle 1, -\alpha_1 \rangle \otimes \cdots \otimes \langle 1, -\alpha_n \rangle$ to the cup product $(\alpha_1) \cdot \cdots \cdot (\alpha_n)$. These group homomorphisms induce an isomorphism of graded rings
\[
	\bigoplus_{n \geq 0} I^n(k)/I^{n+1}(k) \cong \bigoplus_{n \geq 0} H^n(k, \Ztwo).
\]

Now we collect the group homomorphisms $e_n(k)$ into a single map on invariants. The maps $e_n(k) : I^n(k) \to H^n(k, \Ztwo)$ commute with field extensions in the sense that for any extension $k \subseteq k'$, the diagram
\begin{equation*}%\label{eq:Functorialityofen}
	\begin{tikzcd}
	I^n(k) \arrow{r}{e_n(k)} \arrow{d} & H^n(k, \Ztwo) \arrow{d} \\
	I^n(k') \arrow{r}{e_n(k')} & H^n(k', \Ztwo)
	\end{tikzcd}
\end{equation*}
commutes. So for any smooth linear algebraic group $G$, we have a well-defined group homomorphism $e_n : \Inv_{k_0}(G, I^n) \to \Inv_{k_0}(G, H^n)$, and this homomorphism has kernel $\Inv_{k_0}(G, I^{n+1})$. In other words, there is an exact sequence of abelian groups
\[
 0 \to \Inv_{k_0}(G, I^{n+1}) \to \Inv_{k_0}(G, I^{n}) \xrightarrow{e_n} \Inv_{k_0}(G, H^n).
\]

\begin{remark}\label{rmk:en_NonSurj}
It is not clear, in general, if the map $e_n$ on invariants is surjective. 
\end{remark}

One way to interpret Remark \ref{rmk:en_NonSurj} is that the functor $\Inv_{k_0}(G, -)$ is a kind of $\Hom$ functor, and as one would expect, it is left exact but not necessarily right exact. 

\begin{definition}\label{def:lift} 
Let $h \in \Inv_{k_0}(G, H^n)$. A Witt invariant $w \in \Inv_{k_0}(G, I^n)$ is a \emph{lift} of $h$ if $h = e_n \circ w$. 
\end{definition}

For $G$ a Weyl group and $k_0$ a field whose characteristic does not divide $|G|$, we show in Proposition \ref{prop:WeylSurj} that the map $e_n$ on invariants is surjective for all $n$, and we show in Theorem \ref{thm:WeylLiftGen} that $\Inv_{k_0}(G, W)$ is generated as a $W(k_0)$-module by a finite set of lifts of cohomological invariants of $G$. 

\subsection{Trace forms}
\label{subsec:traceform}

For a finite extension of fields $F/k$, the field trace is the map $\Tr_{F/k}: F \to k$ sending $\alpha \in F$ to the trace of the $k$-linear transformation $\beta \mapsto \alpha \beta$. The trace of a field extension $F/k$ determines a symmetric bilinear form $b_{F/k}$ sending $(x, y) \in F \times F$ to $\Tr_{F/k}(xy)$, which is nondegenerate if the extension $F/k$ is separable \cite[Proposition 2.8, p\ 11]{Neukirch1999}. The \emph{trace form} of the extension is the associated quadratic form, sending $x$ to $\Tr_{F/k}(x^2)$. 

A natural setting for the study of the trace form has been algebraic number theory (see \cite{conner1984survey}), where it is connected to another famous invariant:~the discriminant of an extension $K/\mathbb{Q}$ is the determinant of the trace form. However, the trace form makes sense for any finite extension of fields, and this more general perspective lends itself to an algebraic approach. 

To understand the trace form through the lens of Witt invariants, we first take a slight generalization, from separable field extensions to {\' e}tale algebras. Recall that a degree $n$ {\' e}tale algebra over a field $k$ is an $n$-dimensional $k$-algebra $E$ that can be expressed as a product $E \cong k_1 \times \cdots \times k_r$, where $k_1, \ldots, k_r$ are separable field extensions of $k$. By Galois descent we have
\[
	H^1(k, S_n) \leftrightarrow \{ \textup{isomorphism classes of degree } n \textup{ {\' e}tale algebras over } k \}.
\]

For a degree $n$ {\' e}tale algebra $E \cong k_1 \times \cdots \times k_r$ over $k$, the trace form $q_{E/k}$ is the map defined by $(x_1, \ldots, x_r) \mapsto \sum_{i=1}^r \Tr_{k_i/k}(x_i^2)$.
% \[
% 	(x_1, \ldots, x_r) \mapsto \sum_{i=1}^r \Tr_{k_i/k}(x_i^2).
% \]
This definition of the trace form determines a Witt invariant of $S_n$, that is, we have an invariant
\[
	a : H^1(-, S_n) \to W(-)
\]
that sends any degree $n$ {\' e}tale algebra $E$ in $H^1(k, G)$ to the class of the quadratic form $q_{E/k}$ in $W(k)$. For another description of the invariant $a$, using Galois descent, see Section \ref{subsec:Sn_On_traceforms}.

In \cite{Serre03_GMS}, Serre proved that \emph{all} of the Witt invariants of {\' e}tale algebras come from the trace form, in the following sense. For a quadratic form $q$ on an $n$-dimensional vector space $V$ and an integer $0 \leq i \leq n$, the $i$th exterior power $\lambda^i q$ is the quadratic form on $\bigwedge^i V$ corresponding to the symmetric bilinear form 
\[
	(\lambda^d b)(x_1 \wedge \cdots \wedge x_n, y_1 \wedge \cdots \wedge y_n) = \det \Big(b(x_i,y_j)\Big)_{1 \leq i,j \leq n},
\]
see \cite[IX.\S 1.9]{Bourbaki2007}. In particular, $\lambda^0 q = \langle 1 \rangle$, $\lambda^1 q = q$, and for $q = \langle \alpha_1, \ldots, \alpha_n \rangle$ one has
\begin{equation}\label{eq:lambdaformula}
	\lambda^d q = 	\lambda^d(\langle \alpha_1, \ldots, \alpha_n \rangle) = \sum_{\substack{I \subseteq \{1, \ldots, n\} \\ |I| = d}} \langle \prod_{i \in I} \alpha_i \rangle.
\end{equation} 
Serre showed that the Witt invariants of $S_n$ are generated as a free module by the trace form and its exterior powers. 

\begin{theorem}[{\cite[Theorem 29.2, p\ 69]{Serre03_GMS}}]\label{thm:InvOfSn}
$\Inv_{k_0}(S_n, W)$ is a free $W(k_0)$-module with basis $\{a_0, \ldots, a_{\floor{n/2}}\}$, where
\[
	a_i(E) = \lambda^i q_{E/k}
\]
for all fields $k/k_0$ and all {\' e}tale algebras $E \in H^1(k, S_n)$.
\end{theorem}

In Theorems \ref{thm:InvWeylBn}, \ref{thm:InvWeylDn}, and \ref{thm:InvWeylG2}, we generalize Serre's trace form result from {\' e}tale algebras to other classes of algebras arising from Galois cohomology. In particular, we show that a version of Serre's result holds when $S_n$ is replaced by a Weyl group of type $B_n$ (or $C_n$), type $D_n$, or type $G_2$: for these groups, all Witt invariants come from trace forms.

\section{Computations for general \texorpdfstring{$G$}{G}}
\label{sec:Computations}

This section describes several results that can be used to study $\Inv_{k_0}(G, W)$ when $G$ is a smooth linear algebraic group over $k_0$. In Proposition \ref{prop:InvGenMod}, we give sufficient conditions for lifting a set of generators of the cohomological invariants of $G$ to a set of generators for the Witt invariants of $G$. Lemma \ref{lem:ProductInvar} describes the invariants of a product of groups when one of the groups has completely decomposable invariants. In Lemma \ref{lem:StiefelWhitneyLift}, we consider the case $G = O_n$ and observe that the Stiefel-Whitney classes lift to Witt invariants, from which we can conclude that any cohomological invariant that factors through an orthogonal representation lifts to a Witt invariant.

\medskip

As a preliminary we record the following result, which shows that the two possible ways of defining $\Inv_{k_0}(G, H)$ are equivalent. This fact appears to be well-known (for example, it seems to be assumed implicitly in \cite{Serre03_GMS}). Since we did not find a reference in the literature, we record an argument here, which was communicated by Stefan Gille.

\begin{proposition}\label{prop:EquivGradedH}
Let $G$ be a smooth linear algebraic group over $k_0$, and let $H^n$ and $\bigoplus_{n \geq 0} H^n$ be the functors $H^n(-, \Ztwo)$ and $\bigoplus_{n \geq 0} H^n(-, \Ztwo)$ respectively. Then the inclusion map
\[
	\bigoplus_{n \geq 0} \Inv_{k_0}(G, H^n) \to \Inv_{k_0}(G, \bigoplus_{n \geq 0} H^n),
\]
sending $(a_0, a_1, \ldots, a_N, 0, 0, \ldots)$ to $\sum_{n=0}^N a_n$, is an isomorphism.
\end{proposition}

The fact that the inclusion map is an isomorphism does not follow immediately from the definitions:~one could imagine that there might be an invariant $a : H^1(-, G) \to \bigoplus_{n \geq 0} H^n(-, \Ztwo)$ where for each $n$ there exists a field $k/k_0$ and an element $T \in H^1(k, G)$ such that $a(T)$ has degree at least $n$. The content of the proposition is that no such invariant exists.

\begin{proof}[Proof of Proposition \ref{prop:EquivGradedH}]
By \cite[5.3, p\ 12]{Serre03_GMS}, there exists a versal $G$-torsor $P \in H^1(K, G)$, where $K$ is a finitely generated field extension of $k_0$. By \cite[Theorem 12.3, p\ 31]{Serre03_GMS}, two invariants $a, b \in \Inv_{k_0}(G, \bigoplus_{n \geq 0} H^n)$ are equal if and only if $a(P) = b(P)$. 

Now fix $a \in \Inv_{k_0}(G, \bigoplus_{n \geq 0} H^n)$. For $i \geq 0$, let $a_i \in \Inv_{k_0}(G, H^i)$ be the invariant sending $T \in H^1(k, G)$ to the $i$th component of $a(T)$ as an element of $\bigoplus_{n \geq 0} H^n(k, \Ztwo)$. For some $N$, we have in the ring $\bigoplus H^n(K, \Ztwo)$ the equality
\[
	a(P) = \sum_{i=0}^N a_i(P),
\]
which by \cite[Theorem 12.3, p\ 31]{Serre03_GMS} implies that $a$ is equal to the finite sum $\sum_{i=0}^N a_i$.
\end{proof}

\subsection{Lifting cohomological invariants to Witt invariants}
\label{subsec:LiftInvar}

Proposition \ref{prop:InvGenMod} is a kind of ``lifting condition,'' and it is the key step that allows us to translate information about cohomological invariants to information about Witt invariants. In particular, it is essential to the proofs of Theorems \ref{thm:WeylLiftGen}, \ref{thm:InvWeylBn}, and \ref{thm:InvWeylDn}. To state and prove the proposition we first establish some terminology.

For any functor $\FFun : \catFields \to \catAb$, let $\NInv_{k_0}(G, \FFun)$ be the additive group of normalized cohomological invariants, that is, the invariants $a \in \Inv_{k_0}(G, \FFun)$ such that $a(T_k) = 0$ for all $k/k_0$, where $T_k$ is the trivial torsor in $H^1(k, G)$. 
There is a splitting $\Inv_{k_0}(G, \FFun) =  \FFun(k_0) \oplus \NInv_{k_0}(G, \FFun)$ as abelian groups:~each $a \in \Inv_{k_0}(G, \FFun)$ decomposes as $a = a(T_{k_0}) + a'$ where $T_{k_0}$ is the trivial torsor in $H^1(k_0,G)$ and $a'$ is the normalized invariant $a - a(T_{k_0})$. 

We make the following definition, which is a condition that we will use in the proof of Proposition \ref{prop:InvGenMod}. As before we write $H^n(-)$ for the functor $H^n(-, \Ztwo)$, and $H(-)$ for $\bigoplus_{n \geq 0} H^n(-, \Ztwo)$. 

\begin{definition}\label{def:BoundedCohInvar}
Let $G$ be a smooth linear algebraic group and let $n_0$ be an integer. We say that $G$ has $n_0$-\emph{bounded cohomological invariants over} $k_0$ if $\NInv_{k_0}(G, H^n) = 0$ for all $n > n_0$. If there exists an integer $n_0$ such that $G$ has $n_0$-bounded cohomological invariants over $k_0$, we say $G$ has \emph{bounded cohomological invariants over} $k_0$.
\end{definition}

Later, when we consider Weyl groups, we will be able to replace the above condition with an assumption on the characteristic of $k_0$; see Lemma \ref{lem:Weyl_Char_Implies_BddInvar}.

\begin{remark}\label{rmk:FiniteCD_Implies_BddCohInvar}
If $k_0$ has finite $2$-cohomological dimension, then every $G$ has bounded cohomological invariants over $k_0$. This can be seen by evaluating at the versal torsor, similarly to the proof of Proposition \ref{prop:EquivGradedH}: By \cite[5.3, p\ 12]{Serre03_GMS}, there exists a versal $G$-torsor $P$ over a finitely generated extension $K$ of $k_0$. Since $k_0$ has finite $2$-cohomological dimension, so does $K$, and there exists an integer $n_0$ such that $H^n(K, \Ztwo) = 0$ for every $n > n_0$. In particular, if $n > n_0$ and $a \in \NInv_{k_0}(G, H^n)$, then $a(P) = 0$, so $a=0$ by \cite[Theorem 12.3, p\ 31]{Serre03_GMS}.
\end{remark}

\begin{proposition}\label{prop:InvGenMod}
Let $G$ be a smooth linear algebraic group with $n_0$-bounded cohomological invariants over $k_0$. Suppose that for $n = 0, \ldots, n_0$ we have Witt invariants $\{w_{n,i}\}_{\mathcal{I}_n} \subseteq \Inv_{k_0}(G, I^n)$ such that $\bigcup_{n=0}^{n_0} \{e_n \circ w_{n,i}\}_{i \in \mathcal{I}_n}$ generates $\Inv_{k_0}(G, H)$ as an $H(k_0)$-module. Then $\{\langle 1 \rangle\} \cup \bigcup_{n=0}^{n_0} \{w_{n,i}\}_{i \in \mathcal{I}_n}$ generates $\Inv_{k_0}(G, W)$ as a $W(k_0)$-module.
\end{proposition}

\begin{proof} 
Let $u \in \Inv(G, W)$. We will show that $u$ can be expressed as a sum of the form
\begin{equation}\label{eq:ProofOfInvGenModDecomp}
	u = u^{(0)} + u^{(1)} + \cdots + u^{(n_0)} + q,
\end{equation}
where $q \in W(k_0)$ is a constant invariant, and for each $n$, the invariant $u^{(n)} \in \Inv(G, I^n)$ is in the $W(k_0)$-submodule of $\Inv_{k_0}(G, W)$ generated by $\bigcup_{n=0}^{n_0} \{e_n \circ w_{n,i}\}_{i \in \mathcal{I}_n}$.

We begin with an observation about the cohomological invariants of $G$. For $n = 0, \ldots, n_0$ and $i \in \mathcal{I}_n$, write $h_{n,i}$ to denote the cohomological invariant $e_n \circ w_{n,i}$. By assumption, $\Inv_{k_0}(G, H)$ is generated as an $H(k_0)$-module by the set $\bigcup_{n=0}^{n_0} \{h_{n,i}\}_{i \in \mathcal{I}_n}$, so any degree $n$ cohomological invariant of $G$ can be expressed as a sum of the form
\[
	\sum_{m=0}^{n} \sum_{i \in \mathcal{I}_m} c_{n-m,i} h_{m,i}
\]
where $c_{n-m,i}$ is in $H^{n-m}(k_0)$ for every $m$ and every $i$. Any sum of this form can be lifted to a Witt invariant in the following way. Choose quadratic forms $q_{n-m,i} \in I^{n-m}(k_0)$ such that $e_{n-m}(q_{n-m,i}) = c_{n-m,i}$. Then
\[
	e_n \circ \left( \sum_{m=0}^{n} \sum_{i \in \mathcal{I}_m} q_{n-m,i} w_{m,i} \right) =  \sum_{m=0}^{n} \sum_{i \in \mathcal{I}_m} e_{n-m}(q_{n-m,i}) e_{m} \circ w_{m,i} = \sum_{m=0}^{n} \sum_{i \in \mathcal{I}_m} c_{n-m,i} h_{m,i}.
\]

Let $h^{(0)} = e_0 \circ u$. For some coefficients $c_{0,j}^{(0)} \in H^0(k_0)$ and some subset $\mathcal{J}_0^{(0)}$, we have
\[
	h^{(0)} = \sum_{i \in \mathcal{I}_0^{(0)}} c_{0,j}^{(0)} h_{0,j}.
\]
Let $q_{0,j}$ be quadratic forms such that $e_{0,j}(q_{0,j}^{(0)}) = c_{0,j}^{(0)}$ for all $i \in \mathcal{I}_0^{(0)}$, and let $u^{(0)} = \sum_{i \in \mathcal{I}_0^{(0)}} q_{0,j}^{(0)} w_{0,j}$. We have
\[
	e_0 \circ u^{(0)} = \sum_{i \in \mathcal{I}_0^{(0)}} c_{0,j}^{(0)} h_{0,j}.
\]

Setting $w^{(1)} = u - u^{(0)}$, we have $e_0 \circ w^{(1)} = 0$, which implies $w^{(1)} \in \Inv_{k_0}(G, I)$.

We prove the existence of the decomposition (\ref{eq:ProofOfInvGenModDecomp}) by induction on $n$. Assume that there exists $w^{(n)} \in \Inv(G, I^{n})$ such that
\[
 u = u^{(0)} + u^{(1)} + \cdots + u^{(n-1)} + w^{(n)},
\]
where for $0 \leq \ell \leq n-1$ we have
\[
	u^{(\ell)} = \sum_{m=0}^{\ell} \sum_{i \in \mathcal{I}_m} q_{\ell-m, i}^{(\ell)} w_{m, i}.
\]
We have shown that this holds in the $n=0$ case. Setting $h^{(n)} = e_n \circ w^{(n)}$, we have
\[
	h^{(n)} = \sum_{m=0}^{n} \sum_{i \in \mathcal{I}_m} c_{n-m,i}^{(n)} h_{m,i}.
\]
Set 
\[
	u^{(n)} = \sum_{m=0}^{n} \sum_{i \in \mathcal{I}_m} q_{n-m,i}^{(n)} w_{m,i},
\]
where $e_{n-m}(q_{n-m,i}^{(n)}) = c_{n-m,i}^{(n)}$ for all $m$ and all $i$, and set 
\[
	w^{(n+1)} =  w^{(n)} - u^{(n)}.
\]
Since $e_n \circ w^{(n+1)} = 0$, the invariant $w^{(n+1)}$ is in $\Inv(G, I^{n+1})$, and we have $u = u^{(0)} + u^{(1)} + \cdots + u^{(n)} + w^{(n+1)}$.

We have shown that
\[
	u = u^{(0)} + u^{(1)} + \cdots + u^{(n_0)} + w^{(n_0+1)},
\]
where $w^{(n_0+1)} \in \Inv(G, I^{n_0+1})$, and each invariant $u^{(n)}$ is in the $W(k_0)$-submodule of $\Inv_{k_0}(G, W)$ generated by $\bigcup_{n=0}^{n_0} \{w_{n,i}\}_{i \in \mathcal{I}_n}$.

We claim that $w^{(n_0+1)}$ is a constant invariant, and is therefore in the $W(k_0)$-submodule of $\Inv_{k_0}(G, W)$ generated by $\{\langle 1 \rangle\}$. Let $\tilde{w}^{(n_0+1)} = w^{(n_0+1)} - q$, where $q \in W(k_0)$ is the value of $w^{(n_0+1)}$ on the trivial torsor. Since $G$ has $n_0$-bounded cohomological invariants over $k_0$, we have $\NInv(G, H^{n_0+1}) = 0$, and in particular $e_{n_0+1} \circ \tilde{w}^{(n_0+1)} = 0$. This implies that $\tilde{w}^{(n_0+1)}$ is in $\Inv(G, I^{n_0+2})$, and similarly $\tilde{w}^{(n_0+1)}$ is in $\Inv(G, I^{n_0+m})$ for all $m \geq 0$. By the Arason-Pfister Hauptsatz \cite{ArasonPfister1971}, $\bigcap_{n \geq 0} I^n(k) = 0$ for any field $k$. So $\tilde{w}^{(n_0+1)} = 0$, and $w^{(n_0+1)} = q$.
\end{proof}

In Sections \ref{sec:Weyl_CohLifts} and \ref{sec:WInv_TraceForms_BnG2} we apply Proposition \ref{prop:InvGenMod} when $G$ is a Weyl group. In this situation we may simplify our assumption on $k_0$:~the condition of bounded cohomological invariants holds whenever the characteristic of the field does not divide the order of $G$. The following lemma and its proof were communicated by the anonymous referee.

\begin{lemma}\label{lem:Weyl_Char_Implies_BddInvar}
Suppose that $G$ is a Weyl group. If $\textup{char}(k_0)$ does not divide $|G|$, then $G$ has $n_0$-bounded cohomological invariants over $k_0$ for some $n_0$.
\end{lemma}

The proof of the lemma uses Serre's splitting principle. This is \cite[Theorem 25.15, p\ 60]{Serre03_GMS}, but Serre assumes that $k_0$ has characteristic zero. The statement of the splitting principle in \cite[Proposition 2.3, p\ 770]{hirsch2020decomposability} only requires that the characteristic of $k_0$ does not divide $|G|$ and that $H(-)$ is a cycle module; the proof of this version of the splitting principle is given in \cite{GILLE2021}. By assumption $k_0$ has characteristic not equal to $2$, so $H(-)$ is a cycle module. 

\begin{proof}
Let $n_0$ be the largest number such that $G$ contains a subgroup isomorphic to $(\mathbb{Z}/2)^{n_0}$. Suppose $a$ is a normalized cohomological invariant $H^1(-, G) \to H^n(-)$ for some $n > n_0$. For each subgroup $(\mathbb{Z}/2)^{m}$ of $G$, for any $m$, the composition of $H^1(-, (\mathbb{Z}/2)^{m}) \to H^1(-, G)$ with $a$ is an element of $\Inv_{k_0}((\mathbb{Z}/2)^{m}, H^n)$. Using the description of $\Inv_{k_0}((\mathbb{Z}/2)^{m}), H)$ given in \cite[Theorem 16.4, p\ 40]{Serre03_GMS}, we see that the composition must be identically zero, because $a$ is normalized and $m < n$. Since this holds for every abelian subgroup of $G$ generated by reflections, we find by Serre's splitting principle that $a = 0$.
\end{proof}

\subsection{Invariants of a product of groups}
\label{subsec:InvarProducts}

In Section \ref{sec:Weyl_CohLifts} we will consider the Witt invariants of Weyl groups. Recall that any Weyl group $G$ is a product of the form $G_1 \times \cdots \times G_r$, where each $G_i$ is the Weyl group of an irreducible root system. Lemma \ref{lem:ProductInvar} will allow us to extend our results about the Weyl groups of irreducible root systems to all Weyl groups.

We establish some notation that we will use in the statement and proof of Lemma \ref{lem:ProductInvar}. Let $G$ and $G'$ be smooth linear algebraic groups over $k_0$ and $\FFun : \catFields \to \catRings$. We have a product map $\Inv_{k_0}(G, \FFun) \times \Inv_{k_0}(G', \FFun) \to \Inv_{k_0}(G \times G', \FFun)$, defined as follows. Let $a \in \Inv_{k_0}(G, \FFun)$ and $a' \in \Inv_{k_0}(G', \FFun)$. Their product $a \cdot a' \in \Inv_{k_0}(G \times G', \FFun)$ is the invariant sending a pair $(T, T') \in H^1(k, G) \times H^1(k, G') = H^1(k, G \times G')$ to the product $a(T) a'(T) \in \FFun(k)$. 

For any invariant $a \in \Inv_{k_0}(G, \FFun)$ and field extension $k/k_0$, write $\res_{k/k_0}(a) \in \Inv_{k}(G, \FFun)$ for the restriction of $a$ to field extensions of $k$, that is, the invariant defined by $\res_{k/k_0}(a)(T) = a(T)$ for every field extension $K/k$ and torsor $T \in H^1(K, G)$. Following the terminology of \cite{hirsch2020decomposability}, we say that $\Inv_{k_0}(G, \FFun)$ is \emph{completely decomposable} if there is a finite set $\{a_i\}_{i=1}^n$ of $\FFun$-invariants of $G$ such that for every field extension $k/k_0$, the ring of invariants $\Inv_{k}(G, \FFun)$ is a free $\FFun(k)$-module with basis $\left\{\res_{k/k_0}(a_i)\right\}_{i=1}^n$. 

Lemma \ref{lem:ProductInvar} gives a description of the invariants of the product $G \times G'$ in the case that $\Inv_{k_0}(G, \FFun)$ is completely decomposable. 

\begin{lemma}\label{lem:ProductInvar}
Let $G$ and $G'$ be smooth linear algebraic groups over $k_0$ and let $\FFun : \catFields \to \catRings$ be a functor. If $\Inv_{k_0}(G, \FFun)$ is completely decomposable, then the product map
	\begin{align*}
		p : \Inv_{k_0}(G, \FFun) \otimes_{\FFun(k_0)} \Inv_{k_0}(G', \FFun) & \to \Inv_{k_0}(G \times G', \FFun) \\
		a \otimes a' & \mapsto a \cdot a'
	\end{align*}
is an isomorphism.
\end{lemma}

Lemma \ref{lem:ProductInvar} builds directly on earlier results:~see \cite[Exercise 16.5, p\ 40]{Serre03_GMS} and \cite[Proposition 2.5, p\ 771]{hirsch2020decomposability} for special cases of this statement. The following proof was written in collaboration with Saurabh Gosavi. 

\begin{proof}
We will first show that $p$ is surjective. Let $\theta \in \Inv_{k_0}(G \times G', \FFun)$. If $k$ is a field extension of $k_0$ and $T' \in H^1(k, G')$, we have an invariant $\theta(-, T') \in \Inv_{k}(G, \FFun)$ given by the map
\[
	T \mapsto \theta(T, T').
\]
To check that $\theta(-, T')$ is indeed an invariant over $k$ (specifically, that it is a natural transformation), note that if $k \subseteq K \subseteq L$ then we have a commutative diagram
\begin{equation*}
	\begin{tikzcd}
	 H^1(K, G) \times H^1(K, G') \arrow{r}{\theta} \arrow{d} & \FFun(K) \arrow{d} \\
	 H^1(L, G) \times H^1(L, G') \arrow{r}{\theta} & \FFun(L)
	\end{tikzcd}
\end{equation*}
where the vertical map on the left-hand side is the product of the usual maps $H^1(K, G) \to H^1(L, G)$ and $H^1(K, G') \to H^1(L, G')$. Fixing $T' \in H^1(k, G')$ gives the required commutative diagram
\begin{equation*}
	\begin{tikzcd}
	 H^1(K, G) \arrow{r}{\theta(-, T'_K)} \arrow{d} & \FFun(K) \arrow{d} \\
	 H^1(L, G) \arrow{r}{\theta(-, T'_L)} & \FFun(L),
	\end{tikzcd}
\end{equation*}
where $T'_K$ denotes the image of $T'$ under the map $H^1(k, G') \to H^1(K, G')$. 

Since $\Inv_{k_0}(G, \FFun)$ is completely decomposable, there exist invariants $\{a_i\}_{i=1}^n$ such that for every field extension $K$ of $k_0$, $\Inv_{K}(G, \FFun)$ is a free $\FFun(K)$-module with basis $\left\{\res_{K/k_0}(a_{i})\right\}_{i=1}^n$. In particular, for $k$ and $T'$ as above, $\theta(-, T')$ has a unique decomposition of the form
\[
	\theta(-, T') = \sum_{i=1}^{n} b_i(T') \res_{k/k_0}(a_{i})
\]
where $b_1(T'), \ldots, b_n(T')$ are in $\FFun(k)$.

We claim that for each $i$, the assignment $T' \mapsto b_i(T')$ defines an $\FFun$-invariant of $G'$ over $k_0$. Let $T' \in H^1(k, G')$ and let $K$ be a field extension of $k$. We need to check that the image $b_i(T')_K$ of $b_i(T')$ under the map $\FFun(k) \to \FFun(K)$ is equal to $b_i(T'_K)$. 
By restricting invariants to the larger base field $K$, we have in $\Inv_{K}(G, \FFun)$ the equality
\[
	\res_{K/k}(\theta(-,T')) = \res_{K/k}\left( \sum_i b_i(T') \res_{k/k_0}(a_i) \right) = \sum_i b_i(T')_K \res_{K/k_0}(a_i).
\]
On the other hand, by the definition of $\theta(-,T')$ and the definition of the maps $b_i$ respectively, we have
\[
	\res_{K/k} \theta(-,T') = \theta(-,T'_K) = \sum_i b_i(T'_K) \res_{K/k_0} a_i.
\]
Since $\left\{ \res_{K/k_0}(a_i) \right\}_{i=1}^n$ is a basis of $\Inv_K(G, \FFun)$, we have the desired equality $b_i(T')_K = b_i(T'_K)$.

For injectivity, suppose that $p(\sum_\ell \eta_\ell \otimes \eta_\ell') = 0$ for some invariants $\eta_\ell$ of $G$ and $\eta_\ell'$ of $G'$. Write $\eta_\ell = \sum_i b_{\ell,i} a_i$, where $b_{\ell,i} \in \FFun(k_0)$. We have
\[
	\sum_\ell \eta_\ell \otimes \eta_\ell' = \sum_\ell \sum_i b_{\ell,i} a_i \otimes \eta_\ell' = \sum_i a_i \otimes \left( \sum_\ell b_{\ell,i} \eta_\ell' \right).
\]
Applying $p$ gives
\[
	\sum_i a_i \cdot \left( \sum_\ell b_{\ell,i} \eta_\ell' \right) = 0
\]
as an element of $\Inv_{k_0}(G \times G', \FFun)$. For any $k$ and any $T' \in H^1(k, G')$, we get an equation in $\Inv_k(G, \FFun)$ given by
\[
	\sum_i \left( \sum_\ell b_{\ell,i} \eta_\ell'(T') \right) \res_{k/k_0} a_i = 0,
\]
and we know that $\sum_\ell b_{\ell,i} \eta_\ell'(T') \in \FFun(k)$. Since the invariants $\res_{k/k_0}(a_i)$ form a basis for the free $\FFun(k)$-module $\Inv_k(G, \FFun)$, we can conclude that
\[
	\sum_\ell b_{\ell,i} \eta_\ell'(T') = 0
\]
for $i = 0, \ldots, n$. Since this holds for every $T'$ we get
\[
	\sum_i a_i \otimes \left( \sum_\ell b_{\ell,i} \eta_\ell' \right) = \sum_i a_i \otimes 0 = 0,
\]
so $\sum_\ell \eta_\ell \otimes \eta_\ell' = 0$.
\end{proof}

\subsection{Lifting via orthogonal representations}
\label{subsec:LiftOrthReps}

Any orthogonal representation $\rho : G \to O_n$ induces a map $\rho_* : H^1(k, G) \to H^1(k, O_n)$. By composing with an invariant of $O_n$, one obtains an invariant of $G$. As we discuss in Section \ref{sec:Weyl_CohLifts}, it was shown in \cite{hirsch2020decomposability} that almost all of the cohomological invariants of Weyl groups factor through orthogonal representations in this way (as long as $\textup{char}(k_0)$ does not divide $|G|$). In this section we show that $e_d : \Inv_{k_0}(O_n, I^d) \to \Inv_{k_0}(O_n, H^d)$ is surjective for every $d$. As a consequence, any degree $d$ cohomological invariant $h \in \Inv_{k_0}(G, H^d)$ that factors through an orthogonal representation of $G$ can be lifted to a Witt invariant $w \in \Inv_{k_0}(G, I^d)$.

The cohomological invariants and Witt invariants of $O_n$ were computed in \cite{Serre03_GMS}, and we recall them here. Galois descent gives a bijection 
\[
	H^1(k, O_n) \leftrightarrow \left\{\begin{array}{c}
	\textup{isomorphism classes of nondegenerate rank } n \\
	\textup{quadratic forms over } k
	\end{array}\right\}.
\]
By \cite[Theorem 17.3, p\ 41]{Serre03_GMS}, $\Inv_{k_0}(O_n, H)$ is a free $H(k_0)$-module with basis given by the Stiefel-Whitney classes $\{\sw_0, \sw_1, \ldots, \sw_n\}$, where the $d$th Stiefel-Whitney class is defined by
\begin{align*}%\label{eq:DefStiefelWhitney}
	\sw_d : H^1(k, O_n) &\to H^d(k, \Ztwo)\nonumber \\
	\langle \alpha_1, \ldots, \alpha_n \rangle &\mapsto \sum_{1 \leq i_1 < \cdots < i_d \leq n} (\alpha_{i_1}) \cdot \cdots \cdot (\alpha_{i_d}).
\end{align*}
By \cite[Theorem 27.16, p\ 66]{Serre03_GMS}, $\Inv_{k_0}(O_n, W)$ is a free $W(k_0)$-module with basis given by the exterior power invariants $\{\lambda^0, \lambda^1, \ldots, \lambda^n\}$, where $\lambda^i$ is defined, as in (\ref{eq:lambdaformula}), by
\begin{align*}
	\lambda^i : H^1(k, O_n) &\to W(k) \\
	\langle \alpha_1, \ldots, \alpha_n \rangle &\mapsto \sum_{\substack{I \subseteq \{1, \ldots, n\} \\ |I| = d}} \langle \prod_{i \in I} \alpha_i \rangle
\end{align*}
for every field extension $k$ of $k_0$.

Through composition with the maps $e_d$, any Witt invariant of $O_n$ valued in $I^d(-)$ gives a cohomological invariant of $O_n$ valued in $H^d(-, \Ztwo)$. The following lemma goes in the other direction: it describes how to lift the generators of the cohomological invariants of $O_n$ to particular Witt invariants of $O_n$.

\begin{lemma}\label{lem:StiefelWhitneyLift}
The Stiefel-Whitney class $\sw_d : H^1(-, O_n) \to H^d(-, \Ztwo)$ satisfies the equation
\begin{equation*}%\label{eq:StiefelWhitneyLift}
	\sw_d = e_d \circ \left( \sum_{\ell=0}^d (-1)^\ell \binom{n - \ell}{d - \ell} \lambda^\ell \right).
\end{equation*}
In particular, $\sw_d$ lifts to an invariant $H^1(-, O_n) \to I^d(-)$ that is a $W(k_0)$-linear combination of exterior power maps.
\end{lemma}

\begin{proof}
By the definition of $e_d$, we have $\sw_d = e_d \circ t_d$ where $t_d : H^1(-, O_n) \to I^d(-)$ is the map
\[
	\langle \alpha_1, \ldots, \alpha_n \rangle \mapsto \sum_{1 \leq i_1 < \cdots < i_d \leq n} \llangle \alpha_{i_1}, \ldots, \alpha_{i_d} \rrangle.
\]
We will show that $t_d = \sum_{\ell=0}^d (-1)^\ell \binom{n-\ell}{d-\ell} \lambda^\ell$. For any order $d$ subset $I = \{i_1, \ldots, i_d\} \subseteq \{1, \ldots, n\}$, we have
\[
	\llangle \alpha_{i_1}, \ldots, \alpha_{i_d} \rrangle = \prod_{j=1}^d (\langle 1 \rangle - \langle \alpha_{i_j} \rangle) = \sum_{L \subseteq I} (-1)^{|L|} \prod_{i \in L} \langle \alpha_i \rangle.
	%\sum_{\ell=0}^d (-1)^\ell 
\]
% Grouping these terms according to $|L|$ we see that
% \[
% 	\llangle \alpha_{i_1}, \cdots, \alpha_{i_d} \rrangle =
% 		\sum_{\ell=0}^d \sum_{\substack{L \subseteq I \\ |L| = \ell}} (-1)^{\ell} \prod_{i \in L} \langle \alpha_i \rangle =
% 		\sum_{\ell=0}^d (-1)^{\ell} \lambda^\ell \langle \alpha_{i_1}, \cdots, \alpha_{i_d} \rangle.
% \]
So $t_d$ sends $\langle \alpha_1, \ldots, \alpha_n \rangle$ to the sum
\[
	\sum_{\substack{I \subseteq \{1, \ldots, n\} \\ |I| = d}} \sum_{L \subseteq I} (-1)^{|L|} \prod_{i \in L} \langle \alpha_i \rangle.
\]
Let $L$ be an order $\ell$ subset of $\{1, \ldots, n\}$, where $0 \leq \ell \leq d$. The number of order $d$ subsets $I \subseteq \{1, \ldots, n\}$ such that $L \subseteq I$ is $\binom{n-\ell}{d-\ell}$: choosing such an $I$ is equivalent to choosing $d-\ell$ elements of $\{1, \ldots, n\}$ that are not already in $L$. Grouping our terms according to $|L|$ gives
\begin{align*}
	\sum_{\substack{I \subseteq \{1, \ldots, n\} \\ |I| = d}} \sum_{L \subseteq I} (-1)^{|L|} \prod_{i \in L} \langle \alpha_i \rangle
	&=
	\sum_{\ell=0}^d (-1)^\ell \binom{n-\ell}{d-\ell} \sum_{\substack{L \subseteq \{1, \ldots, n\} \\ |L| = \ell}} \prod_{i \in L} \langle \alpha_i \rangle \\
	&= 
	\sum_{\ell=0}^d (-1)^\ell \binom{n-\ell}{d-\ell} \lambda^\ell \langle \alpha_1, \ldots, \alpha_n \rangle.
\end{align*}
This shows that
\begin{equation*}
	t_d = \sum_{\ell=0}^d (-1)^\ell \binom{n-\ell}{d-\ell} \lambda^\ell.
\end{equation*}
\end{proof}

As a consequence of Lemma \ref{lem:StiefelWhitneyLift}, the map
\[
	e_d : \Inv_{k_0}(O_n, I^d) \to \Inv_{k_0}(O_n, H^d)	
\]
is surjective for every $d$. To see this, let $d' = \min(n,d)$ and note that any $h \in \Inv_{k_0}(O_n, H^d)$ is of the form
\[
	h = \sum_{r=0}^{d'} c_{d-r} \sw_r,
\]
for some constant invariants $c_{d-r} \in H^{d-r}(k_0)$. Let $q_{d-r} \in I^{d-r}(k_0)$ be quadratic forms such that $e_{d-r}(q_{d-r}) = c_{d-r}$. By Lemma \ref{lem:StiefelWhitneyLift},
\begin{align*}
	h = \sum_{r=0}^{d'} c_{d-r} \sw_r 
	&= \left( \sum_{r=0}^{d'} e_{d-r}(q_{d-r}) \right) \cdot \left( e_r \circ \sum_{\ell=0}^r (-1)^\ell \binom{n - \ell}{r - \ell} \lambda^\ell \right) \\
	& = e_d \circ \left( \sum_{r=0}^{d'} \sum_{\ell=0}^r q_{d-r} (-1)^\ell \binom{n - \ell}{r - \ell} \lambda^\ell \right),
\end{align*}
so $h$ is in the image of $e_d$. In other words, every degree $d$ cohomological invariant of $O_n$ can be lifted to a Witt invariant.

From the preceding discussion, it follows that any degree $d$ cohomological invariant of $G$ that factors through an orthogonal representation can be lifted to a Witt invariant.

\begin{lemma}\label{lem:LiftInvarFactorOrthRep}
Let $h \in \Inv_{k_0}(G, H)$. If there exists an orthogonal representation $\rho : G \to O_n$ such that $h$ can be written as a composition
\[
	H^1(-, G) \xrightarrow{\rho_*} H^1(-, O_n) \xrightarrow{h'} \bigoplus_{d \geq 0} H^d(-, \Ztwo),
\]
then there exist Witt invariants $w_0, \ldots, w_m$ with $w_d \in \Inv_{k_0}(G, I^d)$ such that $h = \sum_{d=0}^m e_d \circ w_d$. In particular, if $h \in \Inv_{k_0}(G, H^d)$ is a degree $d$ cohomological invariant with $h = h' \circ \rho_*$ for $h' \in \Inv_{k_0}(O_n, H^d)$, then $h$ is in the image of the map $e_d : \Inv_{k_0}(G, I^d) \to \Inv_{k_0}(G, H^d)$.
\end{lemma}

\begin{proof}
Decompose $h$ as $h = h_0 + h_1 + \ldots + h_m$ with $h_d \in \Inv_{k_0}(G, H^d)$ for $0 \leq d \leq m$. Then each $h_d$ factors as $h_d = h'_d \circ \rho_*$. Let $w'_d \in \Inv_{k_0}(O_n, I^d)$ be a lift of $h'_d$. Then $h = \sum_{d=0}^m e_d \circ (w'_d \circ \rho_*)$.
\end{proof}

\section{The Witt invariants of Weyl groups via cohomological invariants}
\label{sec:Weyl_CohLifts}

Let $G$ be a Weyl group and $k_0$ a field of characteristic not dividing $|G|$. Hirsch showed that $\Inv_{k_0}(G, H)$ is completely decomposable \cite[Corollary 2.6, p\ 772]{hirsch2020decomposability}, and computed explicit generators for $\Inv_{k_0}(G, H)$ as an $H(k_0)$-module.\footnote{Another description of the $H(k_0)$-module $\Inv_{k_0}(G, H)$, which only assumes that $\textup{char}(k_0)$ is not equal to $2$, can be found in \cite{serre2018cohomological}.} The main result of this section is Theorem \ref{thm:WeylLiftGen}, which says that the $W(k_0)$-module $\Inv_{k_0}(G, W)$ is generated by a finite set of lifts of cohomological invariants. As an intermediate step we prove Proposition \ref{prop:WeylSurj}, which says that the map $e_n$ on invariants of $G$ is surjective. Deducing Theorem \ref{thm:WeylLiftGen} from Proposition \ref{prop:WeylSurj} is an application of the lifting condition described in Proposition \ref{prop:InvGenMod}. 

Most of the generators described in \cite{hirsch2020decomposability} come from composing orthogonal representations $G \to O_n$ with Stiefel-Whitney classes, as Hirsch notes in his introduction \cite[p\ 766]{hirsch2020decomposability}. For the particular descriptions of the generators for Weyl groups of irreducible root systems, see \cite[Proposition 4.1, p\ 780]{hirsch2020decomposability} for type $A_n$ (originally due to Serre {\cite[Theorem 29.2, p\ 69]{Serre03_GMS}}), \cite[Corollary 5.7, p\ 790]{hirsch2020decomposability} for type $B_n$, \cite[Corollary 7.5, p\ 799]{hirsch2020decomposability} for type $D_n$, \cite[Section 8, p\ 800-805]{hirsch2020decomposability} for types $E_6$, $E_7$, and $E_8$, \cite[Proposition 6.1, p\ 793]{hirsch2020decomposability} for type $F_4$, and \cite[Section 3.3, p\ 774-775]{hirsch2020decomposability} for type $G_2$. The only generators that do not factor through Stiefel-Whitney classes are particular invariants for types $D_{2n}$, $E_7$, and $E_8$, and these invariants factor through powers of the fundamental ideal $I^d(-)$: see \cite[p\ 795-796]{hirsch2020decomposability}, \cite[Lemma 8.1, p\ 802]{hirsch2020decomposability}, and \cite[Lemma 8.3, p\ 805]{hirsch2020decomposability} for descriptions of these invariants.

\begin{proposition}\label{prop:WeylSurj} 
Let $G$ be a Weyl group and suppose that $\textup{char}(k_0)$ does not divide $|G|$. Then the map $e_n : \Inv_{k_0}(G, I^n) \to \Inv_{k_0}(G, H^n)$ is surjective for all $n$. 
\end{proposition}

\begin{proof}
Any Weyl group $G$ decomposes as a product $G_1 \times \cdots \times G_r$, where $G_1, \ldots, G_r$ are Weyl groups of irreducible root systems. By Lemma \ref{lem:ProductInvar}, 
\[
	\Inv_{k_0}(G, H) \cong \Inv_{k_0}(G_1, H) \otimes_{H(k_0)} \cdots \otimes_{H(k_0)} \Inv_{k_0}(G_r, H).
\]
Note that since $\Inv_{k_0}(G_i, H)$ is completely decomposable by \cite[Theorem 2.2, p\ 769]{hirsch2020decomposability}, $\Inv_{k_0}(G, H)$ is a free $H(k_0)$-module generated by the products of the generators of the modules $\Inv_{k_0}(G_i, H)$. 

By the classification of cohomological invariants given in \cite{hirsch2020decomposability} (see the discussion preceding Proposition \ref{prop:WeylSurj}), if $G$ is a Weyl group and $\textup{char}(k_0)$ does not divide the order of $G$, then $\Inv_{k_0}(G, H)$ is generated as an $H(k_0)$-algebra by at most two types of invariants:~those that factor through an orthogonal representation $G \to O_n$, and those that can be described by composing $e_n$ with a Witt invariant. By Lemma \ref{lem:LiftInvarFactorOrthRep}, any $h \in \Inv_{k_0}(G, H)$ that factors through an orthogonal representation can be expressed as a sum $h = \sum_{d=0}^m e_d \circ w_d$ for some Witt invariants $w_d \in \Inv_{k_0}(G, I^d)$. Finally, note that if $a_n \in \Inv_{k_0}(G, H^n)$ and $b_m \in \Inv_{k_0}(G, H^m)$, then $e_{n+m}(a_n b_m) = e_n(a_n) e_m(b_m)$: in other words, if two Witt invariants lift to cohomological invariants, then the lift of their product is the product of their lifts. 
\end{proof}

In particular, by Lemma \ref{lem:Weyl_Char_Implies_BddInvar} and Proposition \ref{prop:InvGenMod}, for any Weyl group $G$ and field $k_0$ whose characteristic does not divide the order of $G$, we can describe the Witt invariants of $G$ by lifting cohomological invariants. 

\begin{theorem}\label{thm:WeylLiftGen}
Let $G$ be a Weyl group and suppose that $\textup{char}(k_0)$ does not divide $|G|$. Then $\Inv_{k_0}(G, W)$ is generated as a $W(k_0)$-module by a finite set of lifts of cohomological invariants of $G$.
\end{theorem}

\section{Witt invariants of Weyl groups in terms of trace forms}
\label{sec:WInv_TraceForms_BnG2}

In this section we give detailed descriptions of the Witt invariants of particular Weyl groups in terms of trace forms and their exterior powers, proving an analogue of Serre's trace form result (Theorem \ref{thm:InvOfSn}). In particular, we show in Theorems \ref{thm:InvWeylBn} and \ref{thm:InvWeylDn} that if $G$ is a Weyl group of type $B_n$ or type $D_n$, and $k_0$ is a field whose characteristic is zero or greater than $n$, then the Witt invariants of $G$ over $k_0$ are generated by trace forms and their exterior powers. We also show, in Theorem \ref{thm:InvWeylG2}, that the Witt invariants of the Weyl group of type $G_2$ are generated by trace forms, over any field $k_0$ of characteristic not equal to $2$. Weyl groups of other types are discussed in Remark \ref{rem:OtherExceptionalTypes}.

\subsection{Trace forms via the inclusion \texorpdfstring{$S_n \hookrightarrow O_n$}{Sn to On}}
\label{subsec:Sn_On_traceforms}

As an invariant, the trace form has a simple interpretation in terms of Galois descent. For details on Galois descent see  \cite[\S 2.3]{gillesamuely2017central}, \cite[III.\S 1]{Serre1997} or \cite[\S 18]{KMRT1998}. 

Galois descent gives identifications
\[
	H^1(k, S_n) \leftrightarrow \{ \textup{isomorphism classes of degree } n \textup{ {\' e}tale algebras over } k \}
\]
and
\[
	H^1(k, O_n) \leftrightarrow \left\{\begin{array}{c}
	\textup{isomorphism classes of nondegenerate rank } n \\
	\textup{quadratic forms over } k
	\end{array}\right\}.
\]

We review some of the details of these two identifications, which will also serve to fix our notation for the proof of Lemma \ref{lem:Sn_to_On_is_trace}. Let $E_0 = k \times \cdots \times k$ be the split degree $n$ {\' e}tale $k$-algebra. The group $S_n$ acts on $E_0$ by permuting components. For the purpose of Galois descent, the important fact is that $S_n$ is isomorphic to the automorphism group of $E_0 \otimes_k k^{\sep}$. Let $V_0$ be the $n$-dimensional vector space $k \oplus \cdots \oplus k$ equipped with the quadratic form $q_0 = \langle 1, \ldots, 1 \rangle$, with the orthogonal group $O_n$ acting in the usual way. The set $H^1(k, S_n)$ is identified with the set of twisted forms of $E_0$, and $H^1(k, O_n)$ with the set of twisted forms of $(V_0, q_0)$.

Let $\iota_n : S_n \to O_n$ be the orthogonal representation in which $S_n$ acts on $k^{\oplus n}$ by permuting components. For example, $\iota_3 : S_3 \to O_3$ sends the permutation $(123)$ to the matrix
\[
	\begin{pmatrix}
	0 & 0 & 1 \\
	1 & 0 & 0 \\
	0 & 1 & 0 \\
	\end{pmatrix}.
\]

% The group $S_n$ acts on the vector space $V = k \eps_1 \oplus \cdots \oplus k \eps_n$ by permuting the basis elements $\{\eps_i\}_{i=1}^n$. This action defines an orthogonal representation $\iota_n : S_n \to O_n$. 
Recall that if $G$ and $G'$ are smooth linear algebraic groups over $k$ and $\varphi : G \to G'$ is a $\Gamma_k$-equivariant group homomorphism, then there is an induced map $\varphi_* : H^1(k, G) \to H^1(k, G')$ sending the class of the cocycle $c : \Gamma_k \to G$ to the class of the cocycle $\varphi \circ c : \Gamma_k \to G'$.

\begin{lemma}\label{lem:Sn_to_On_is_trace}
The map
\[
	(\iota_{n})_* : H^1(k, S_n) \to H^1(k, O_n)
\]
sends each degree $n$ {\' e}tale algebra $E$ to its trace form $q_{E/k}$.
\end{lemma}

\begin{proof}
Let $c : \Gamma_k \to S_n$ be a cocycle. Let $E$ be degree $n$ {\' e}tale algebra corresponding to $c$ and let $(V, q)$ be the quadratic space corresponding to $(\iota_{n})_*(c)$. As described above, let $E_0$ be the split degree $n$ {\' e}tale algebra, let $V_0$ be the vector space $k^{\oplus n}$, and let $q_0$ be the quadratic form $\langle 1, \ldots, 1 \rangle$ on $V_0$.

Under the map $\iota_n$, the group $S_n$ acts on $V_0 = k^{\oplus n}$ by permuting the basis, so $V_0 \otimes_k k^{\sep} = E_0 \otimes_k k^{\sep}$ as vector spaces with $S_n$-action. In particular, we can replace $(V_0, q_0)$ with $(E_0, q_0)$, and $q$ is the quadratic form on $E$ that one obtains by twisting $(E_0, q_0)$ by the cocycle $c$.

The trace form on $E_0 \otimes_k k^s$ is $\langle 1, \ldots, 1 \rangle = q_0 \otimes_k k^s$. For any $x \in E$ we have that
\[
	q(x) = (q_0 \otimes_k k^s)(x) = \Tr_{E_0 \otimes_k k^s/k^s}(x^2) = \Tr_{E \otimes_k k^s/k^s}(x^2) = \Tr_{E/k}(x^2).
\]
This shows that $q$ is indeed the trace form on $E$. 
\end{proof}

Lemma \ref{lem:Sn_to_On_is_trace} says that the trace form is the descended version of the canonical embedding of the symmetric group into the orthogonal group. This adds another layer of interpretation to our discussion of trace forms in the context of Witt invariants. For example, in Theorem \ref{thm:InvOfSn}, Serre showed that $\Inv_{k_0}(S_n, W)$ is generated by the trace form invariant $a_1 : H^1(-, S_n) \to W(-)$ and its exterior powers. 
%, where $a_1$ sends each degree $n$ {\' e}tale algebra $E \in H^1(k, S_n)$ to the class of its trace form $q_{E/k}$ in $W(k)$. 
Lemma \ref{lem:Sn_to_On_is_trace} says that $a_1$ is exactly the invariant
\[
	H^1(-, S_n) \xrightarrow{(\iota_{n})_*} H^1(-, O_n) \to W(-),
\]
where the map $H^1(-, O_n) \to W(-)$ sends each nondegenerate degree $n$ quadratic form $q \in H^1(k, O_n)$ to its class in $W(k)$.

If $G$ is any finite group, then by Lemma \ref{lem:Sn_to_On_is_trace} any group homomorphism $\varphi : G \to S_n$ determines a trace form invariant of $G$. Explicitly, the map $\varphi : G \to S_n$ induces a map from $G$-torsors to {\' e}tale algebras $\varphi_* : H^1(-, G) \to H^1(-, S_n)$. Composing with $(\iota_{n})_*$ gives a Witt invariant
\[
	H^1(-, G) \xrightarrow{\varphi_*} H^1(-, S_n) \xrightarrow{(\iota_{n})_*} H^1(-, O_n) \to W(-).
\]
This composition is a trace form invariant, in the sense that it is an invariant that assigns the trace form of an {\' e}tale algebra to each $G$-torsor. We will discuss some specific invariants of this form in the proof of Theorem \ref{thm:InvWeylBn}.

\subsection{Witt invariants of \texorpdfstring{$\Weyl(G_2)$}{Weyl(G2)}}
\label{subsec:winv_G2}

As a consequence of Lemma \ref{lem:ProductInvar} and Theorem \ref{thm:InvOfSn}, we have the following description for the Witt invariants of $\Weyl(G_2) \cong S_2 \times S_3$. Note that this result does not require any assumptions on $k_0$ other than our usual assumption that the characteristic is not equal to $2$.

\begin{theorem}\label{thm:InvWeylG2} 
Let $a_1^{(2)}$ be the trace form invariant in $\Inv_{k_0}(S_2, W)$ and $a_1^{(3)}$ be the trace form invariant in $\Inv_{k_0}(S_3, W)$. The $W(k_0)$-module of invariants $\Inv_{k_0}(\Weyl(G_2), W)$ is a free module with basis $\left\{\langle 1 \rangle, a_1^{(2)}, a_1^{(3)}, a_1^{(2)} \cdot a_1^{(3)}\right\}$. 
\end{theorem}

In particular, $\Inv_{k_0}(\mathcal{W}(G_2), W)$ is generated as a $W(k_0)$-algebra by trace forms. 

\subsection{Witt invariants of \texorpdfstring{$\Weyl(B_n) \cong \Weyl(C_n)$}{Weyl(Bn) = Weyl(Cn)}}
\label{subsec:winv_Bn}

Let $\Weyl(B_n)$ denote the Weyl group of type $B_n$ (which is isomorphic to the Weyl group of type $C_n$). 
As a group, $\Weyl(B_n)$ is isomorphic to the wreath product $(\Ztwo) \wr S_n$, that is, the semidirect product $(\Ztwo)^n \rtimes S_n$ where the symmetric group acts on $(\Ztwo)^n$ by permuting components. Write $s_i$ for the generator of the $i$th copy of $\Ztwo$ in $(\Ztwo) \wr S_n$. An element of the wreath product is of the form
\[
	\sigma \cdot \prod_{i \in I} s_i,
\]
where $\sigma \in S_n$, $I \subseteq \{1, \ldots, n\}$, and multiplication in the group is defined by
\[
	\left(\sigma \cdot \prod_{i \in I} s_i \right) \cdot \left(\tau \cdot \prod_{j \in J} s_{j} \right) = \sigma \tau \cdot \prod_{i \in I} s_{\tau^{-1}(i)} \prod_{j \in J} s_{j} .
\]
For later use (in Theorem \ref{thm:InvWeylBn}), we note that since the order of $\Weyl(B_n)$ is $2^n \cdot n!$, the characteristic of $k_0$ does not divide $|\Weyl(B_n)|$ if and only if the characteristic is zero or greater than $n$.

If $K = k_1 \times \cdots \times k_r$ is an {\' e}tale $k$-algebra, where $k_i$ is a separable field extension of $k$ for each $i$, then an {\' e}tale $k$-algebra $L \supseteq K$ is said to be a \emph{quadratic} extension of $K$ if it is of constant rank $2$, that is, if $L = \ell_1 \times \cdots \times \ell_r$ where $\ell_i \supseteq k_i$ is a degree $2$ {\' e}tale $k_i$-algebra for each $i$. 

\begin{lemma}\label{lem:LemmaH1Bn}
There is a natural bijection
\begin{equation}\label{eq:TypeBnDescription}
	H^1(k, (\Ztwo) \wr S_n) \leftrightarrow 
	\left\{ 
		\begin{array}{c}
		\textup{isomorphism classes of pairs } (K, L) \\ \textup{ of {\' e}tale } k\textup{-algebras } \textup{with } L \supseteq K, \\ \dim_k K = n, \ \textup{and } L/K \textup{ quadratic}
	\end{array} 
	\right\}.
\end{equation}
\end{lemma}

\begin{proof}

We will prove this using Galois descent. As a base point, consider the pair
\[
	K_0 = k \eps_1 \times \cdots \times k \eps_n 
\]
and
\[
	L_0 = k \eps_1^+ \times k \eps_1^- \times \cdots \times k \eps_n^+ \times k \eps_n^-,
\]
where the $\eps_i$ (respectively $\eps_i^\pm$) are primitive orthogonal idempotents, and $K_0 \hookrightarrow L_0$ via $\eps_i \mapsto \eps_i^+ + \eps_i^-$. The pair $(K_0, L_0)$ is an element of the right-hand side in (\ref{eq:TypeBnDescription}). 

Now we consider the automorphisms of the pair $(K_0 \otimes_k k^{\sep}, L_0 \otimes_k k^{\sep})$ over $k^{\sep}$. An automorphism of this pair is, by definition, an automorphism of $L_0 \otimes_k k^{\sep}$ that restricts to an automorphism of $K_0 \otimes_k k^{\sep}$. Such an automorphism corresponds to a permutation of the $2n$ idempotents $\eps_1^+, \eps_1^-, \ldots, \eps_n^+, \eps_n^-$ that also permutes the $n$ idempotents $\eps_1, \ldots, \eps_n$. The group $(\Ztwo) \wr S_n$ is exactly the group of automorphisms of the set $\{\eps_1^+, \eps_1^-, \ldots, \eps_n^+, \eps_n^-\}$ that send each unordered pair of the form $\{\eps_i^+, \eps_i^-\}$ to another pair $\{\eps_j^+, \eps_j^-\}$. So, by Galois descent,
\[
	H^1(k, (\Ztwo) \wr S_n) \leftrightarrow 
	\left\{ \textup{twisted forms of } (K_0, L_0) \right\}.
\]

Taking an element of the right-hand side of (\ref{eq:TypeBnDescription}) and extending scalars to $k^{\sep}$ gives a pair isomorphic to $(K \otimes_k k^{\sep}, L_0 \otimes_k k^{\sep})$. For the other direction, let $(K, L)$ be a twisted form of $(K_0, L_0)$. The fact that $L$ and $K$ are {\' e}tale $k$-algebras follows from the isomorphisms $L \otimes_k k^{\sep} \cong L_0 \otimes_k k^{\sep} \cong (k^{\sep})^{2n}$ and $K \otimes_k k^{\sep} \cong K_0 \otimes_k k^{\sep} \cong (k^{\sep})^{n}$. Since the extension $L \otimes_k k^{\sep} \supseteq K \otimes_k k^{\sep}$ is of constant rank $2$, so is $L \supseteq K$.
\end{proof}

For a pair $(K, L)$ in $H^1(k, (\Ztwo) \wr S_n)$, where $K$ is a degree $n$ {\' e}tale $k$-algebra and $(K, L)$ is a quadratic extension, we have a rank $2n$ trace form $q_{L/k}: x \mapsto \Tr_{L/k}(x^2)$ and a rank $n$ trace form $q_{K/k} : x \mapsto \Tr_{K/k}(x^2)$. Let $a_{L/k}$ and $a_{K/k}$ be the Witt invariants sending $(K, L)$ to $q_{L/k}$ and $q_{K/k}$ respectively.

\begin{theorem}\label{thm:InvWeylBn} 
Suppose that $\textup{char}(k_0)$ does not divide $|\Weyl(B_n)|$, that is, $\textup{char}(k_0)$ is zero or greater than $n$. Then $\Inv_{k_0}(\Weyl(B_n), W)$ is generated as a $W(k_0)$-algebra by $a_{L/k}$, $a_{K/k}$, and their exterior powers.
\end{theorem}

In other words, over a field $k_0$ of appropriate characteristic, the $W(k_0)$-algebra of Witt invariants of $\Weyl(B_n)$ is generated by trace form invariants and their exterior powers.

To prove Theorem \ref{thm:InvWeylBn}, we will lift the cohomological invariants of $\Weyl(B_n)$ computed by Hirsch \cite[Corollary 5.7, p\ 790]{hirsch2020decomposability}, then use Lemma \ref{lem:Weyl_Char_Implies_BddInvar} and Proposition \ref{prop:InvGenMod}. Before we begin the proof we describe these cohomological invariants. For $0 \leq d \leq n$, Hirsch defines the cohomological invariant $u_d$ to be the composition
\[
	u_d : H^1(-, \Weyl(B_n)) \xrightarrow{\rho_*} H^1(-, S_n) \xrightarrow{(\iota_{n})_*} H^1(-, O_n) \xrightarrow{\swmod_d} H^d(-, \Ztwo),
\]
where $\rho_*$ is the map on cohomology induced by the group homomorphism $\rho : (\Ztwo) \wr S_n \to S_n$ sending $\sigma \cdot \prod_{i \in I} s_{i}$ to $\sigma$, the map $(\iota_{n})_*$ is induced by the usual representation $\iota_{n} : S_n \to O_n$ in which $S_n$ acts by permuting the standard basis for $k^n$, and $\swmod_d$ is the modified Stiefel-Whitney class
\[
	\swmod_d(q) = \begin{cases}
		\sw_d(q) & \textup{ if } d \textup{ is odd} \\
		\sw_d(q) + (2) \cdot w_{d-1}(q) & \textup{ if } d \textup{ is even}.
		\end{cases}
\]
Hirsch also defines invariants $v_d'$ via the group homomorphism $\rho_2 : (\Ztwo) \wr S_n \to S_{2n}$ sending 
\[
	\sigma \cdot \prod_{i \in I} s_{i} \mapsto \sigma \cdot (\sigma + n) \cdot \prod_{i \in I} (i, i+n),
\]
where 
\[
	(\sigma + n): i \mapsto 
	\begin{cases}
	i & \textup{ if } i \leq n, \\
	n + \sigma(i-n) & \textup{ if } n < i \leq 2n.
	\end{cases}
\]
The invariant $v_d'$ is the composition
\[
	v_d' : H^1(-, \Weyl(B_n)) \xrightarrow{(\rho_2)_*} H^1(-, S_{2n}) \xrightarrow{(\iota_{2n})_*} H^1(-, O_{2n}) \xrightarrow{\swmod_d} H^d(-, \Ztwo)
\]

Finally, $v_d$ is defined inductively by $v_d = v_d' + \sum_{0 \leq i \leq d-1} u_{d-i} v_i$. Hirsch gives the following decomposition for the cohomological invariants of the Weyl group of type $B_n$.
\begin{theorem}[{\cite[Corollary 5.7, p\ 790]{hirsch2020decomposability}}]\label{thm:HirschBn}
Let $k_0$ be a field whose characteristic does not divide $|\Weyl(B_n)|$. Then $\Inv_{k_0}(\Weyl(B_n), H)$ is a free $H(k_0)$-module with basis 
\[
	\{u_{d-r}v_r : \max(0, 2d-n) \leq r \leq d \leq n \}.
\]
\end{theorem}

\begin{proof}[Proof of Theorem \ref{thm:InvWeylBn}]

Applying Lemma \ref{lem:Weyl_Char_Implies_BddInvar} and Proposition \ref{prop:InvGenMod} to the generators described in Theorem \ref{thm:HirschBn}, we see that it is enough to show that for each $d$ and each $r$, we can lift $u_{d-r}$ and $v_r$ to elements of the $W(k_0)$-algebra generated by $a_{L/k}$, $a_{K/k}$, and their exterior powers. Each $v_r$ is in the algebra generated by the invariants $u_d$ and $v_d'$ for $0 \leq d \leq n$, and we will show that for all $0 \leq d \leq n$, the invariants $u_d$ and $v_d'$ have lifts in the $W(k_0)$-algebra generated by $a_{L/k}$, $a_{K/k}$, and their exterior powers. The argument has three parts.
\begin{enumerate}[(1)]
	\item\label{item:iota_trace_K} $(\iota_{n})_* \circ \rho_*$ sends each pair $(K, L)$ in $H^1(k, \Weyl(B_n))$ to the trace form $q_{K/k}$.
	\item\label{item:iota_trace_L}  $(\iota_{2n})_* \circ (\rho_2)_*$ sends each pair $(K, L)$ in $H^1(k, \Weyl(B_n))$ to the trace form $q_{L/k}$.
	\item\label{item:lift_modsw}  For any $d \leq m$, there exist $q_{m,1}, \ldots, q_{m,d} \in W(k_0)$ such that $\swmod_d : H^1(-, O_m) \to H^d(-, \Ztwo)$ can be expressed as a sum of the form
	\[
		\swmod_d = e_d \circ \left( \sum_{\ell=0}^d q_{m,\ell} \lambda^\ell \right).
	\] 
\end{enumerate}

Combining these three statements, we have that
\[
	u_d = \swmod_d \circ (\iota_{n})_* \circ \rho_* = e_d \circ \left(\sum_{\ell=0}^d q_{m,\ell} \cdot \lambda^\ell a_{K/k} \right)
\]
and
\[
	v_d' = \swmod_d \circ (\iota_{2n})_* \circ (\rho_2)_* = e_d \circ \left(\sum_{\ell=0}^d q_{m,\ell} \cdot \lambda^\ell a_{L/k} \right),
\]
so for every $0 \leq d \leq n$, the invariants $u_d$ and $v_d'$ are in the algebra generated by $a_{K/k}$, $a_{L/k}$, and their exterior powers. 

To prove parts (\ref{item:iota_trace_K}) and (\ref{item:iota_trace_L}), we first note that by Lemma \ref{lem:Sn_to_On_is_trace}, for any $m$ the map $(\iota_m)_*$ sends each degree $m$ {\' e}tale algebra $E \in H^1(k, S_m)$ to its trace form $q_{E/k} \in H^1(k, O_m)$. It remains to show that $\rho_*$ sends $(K, L)$ to $K$ and $(\rho_2)_*$ sends $(K, L)$ to $L$. 

Recall that $\rho : (\Ztwo) \wr S_n \to S_n$ is the map $\sigma \cdot \prod_{i \in I} s_{i} \mapsto \sigma$. By Lemma \ref{lem:LemmaH1Bn}, $H^1(k, (\Ztwo) \wr S_n)$ is identified with the set of twisted forms of the pair of algebras
\[
	k \eps_1^+ \times k \eps_1^- \times \cdots \times k \eps_n^+ \times k \eps_n^- \supseteq k \eps_1 \times \cdots \times k \eps_n
\]
where $\eps_i = \eps_i^+ + \eps_i^-$ for all $i$. At the separable closure the group $(\Ztwo) \wr S_n$ acts by
\[
	\sigma(\eps_i^+) = \eps_{\sigma(i)}^+, \quad \sigma(\eps_i^-) = \eps_{\sigma(i)}^-
\]
and 
\[
	s_j(\eps_i^+) =
	\begin{cases}
		\eps_i^+ & \textup{ if } i \neq j \\
		\eps_i^- & \textup{ if } i = j,
	\end{cases} \
	\qquad
	s_j(\eps_i^-) =
	\begin{cases}
		\eps_i^- & \textup{ if } i \neq j \\
		\eps_i^+ & \textup{ if } i = j.
	\end{cases} \
\]

Given a cocycle class $c \in H^1(k, (\Ztwo) \wr S_n)$ corresponding to a pair $(K, L)$, the degree $n$ {\'e}tale algebra $K$ is obtained by twisting $k \eps_1 \times \cdots \times k \eps_n$ by $c$. For each $\gamma \in \Gamma_k$, the element $c_\gamma = \sigma \cdot \prod_I s_i \in (\Ztwo) \wr S_n$ acts on $k^{\sep} \eps_1 \times \cdots \times k^{\sep} \eps_n$ by sending each $\eps_i$ to $\eps_{\sigma(i)}$. This is exactly the action of $\rho(c_\gamma) \in S_n$ on the split degree $n$ {\' e}tale algebra. So the cocycle $\rho_*(c)$ corresponds to the algebra $K$ in $H^1(k, S_n)$, that is, $\rho_*$ sends $(K, L)$ to $K$.

A similar argument shows that $(\rho_2)_*$ sends $(K, L)$ to $L$. By relabeling idempotents, we identify $k \eps_1^+ \times k \eps_1^- \times \cdots \times k \eps_n^+ \times k \eps_n^-$ with the degree $2n$ split {\' e}tale algebra $k f_1 \times k f_2 \times \cdots \times k f_{2n}$ by $\eps_i^+ \mapsto f_i$ and $\eps_i^- \mapsto f_{n+i}$. At the separable closure, $\sigma \cdot \prod_I s_i \in (\Ztwo) \wr S_n$ acts on the idempotents of $\prod_{j=1}^{2n} k^{\sep} f_j$ by the permutation $\rho_2(\sigma \cdot \prod_I s_i) \in S_{2n}$. In particular, if the cocycle $c$ corresponds to $(K, L)$ in $H^1(k, (\Ztwo) \wr S_n)$, then $(\rho_2)_*(c)$ corresponds to $L$ in $H^1(k, S_{2n})$.

Part (\ref{item:lift_modsw}) follows almost immediately from Lemma \ref{lem:StiefelWhitneyLift}. For the modified Stiefel-Whitney class, we note that $e_1 \circ \llangle 2 \rrangle = (2)$. We have $e_i(q_i) e_j(q_j) = e_{i+j}(q_i \otimes q_j)$ for any $q_i \in I^i(k)$ and $q_j \in I^j(k)$. In particular, the invariant
\[
	\widetilde{t_d} = 
	\begin{cases}
	\sum_{\ell=0}^d (-1)^\ell \binom{n-\ell}{d-\ell} \lambda^\ell & \textup{ if } d \textup{ is odd} \\
	\sum_{\ell=0}^d (-1)^\ell \binom{n-\ell}{d-\ell} \lambda^\ell + 
	\sum_{\ell=0}^{d-1} (-1)^\ell \binom{n-\ell}{d-1-\ell} \llangle 2 \rrangle \lambda^\ell & \textup{ if } d \textup{ is even}
	\end{cases}
\]
satisfies $e_d \circ \widetilde{t_d} = \swmod_d$.
\end{proof}

\subsection{Witt invariants of \texorpdfstring{$\Weyl(D_n)$}{Weyl(Dn)}}
\label{subsec:winv_Dn}

Let $\Weyl(D_n)$ be the Weyl group of type $D_n$, with $n \geq 2$. This group is an index $2$ subgroup of $\Weyl(B_n)$, consisting of the elements $\sigma \cdot \prod_{i \in I} s_i \in \Weyl(B_n)$ such that $|I|$ is even. The inclusion map $\varphi : \Weyl(D_n) \to \Weyl(B_n)$ induces a map on cohomology $\varphi_* : H^1(k, \Weyl(D_n)) \to H^1(k, \Weyl(B_n))$. Note (for use in Theorem \ref{thm:InvWeylDn}) that $|\Weyl(D_n)| = 2^{n-1} \cdot n!$, so the characteristic of $k_0$ does not divide $|\Weyl(D_n)|$ if and only if the characteristic is zero or greater than $n$.

Recall from Section \ref{subsec:winv_Bn} that $\Inv_{k_0}(\Weyl(B_n), W)$ is generated as a $W(k_0)$-algebra by two trace form invariants, $a_{L/k}$ and $a_{K/k}$, and their exterior powers. By composing with $\varphi_*$, one has trace form invariants in $\Inv_{k_0}(\Weyl(D_n), W)$ given by 
\[
	a'_{L/k} = a_{L/k} \circ \varphi_* \quad \textup{and} \quad a'_{K/k} = a_{K/k} \circ \varphi_*.
\] 
For $n$ even, we define an additional invariant $r \in \Inv_{k_0}(\Weyl(D_n), W)$, following \cite[p\ 795]{hirsch2020decomposability}. The group $\Weyl(D_n)$ decomposes as a semidirect product $\Weyl(D_n) \cong (\Ztwo)^{n-1} \rtimes S_n$, so the set of cosets $\Weyl(D_n)/S_n$ has order $2^{n-1}$. The left action of $\Weyl(D_n)$ on these cosets induces a map $\Weyl(D_n) \to S_{2^{n-1}}$. Composing with $\iota_{2^{n-1}}$, we have a Witt invariant
\[
	r : H^1(-, \Weyl(D_n)) \to H^1(-, S_{2^{n-1}}) \xrightarrow{(\iota_{2^{n-1}})_*} H^1(-, O_{2^{n-1}}) \to W(-).
\]
Note that $r$ is a trace form invariant: an element $T \in H^1(k, \Weyl(D_n))$ is first sent to a degree $2^{n-1}$ {\' e}tale algebra $E_T \in H^1(k, S_{2^{n-1}})$, and by Lemma \ref{lem:Sn_to_On_is_trace}, the map $(\iota_{2^{n-1}})_*$ sends $E_T$ to its trace form $q_{E_T/k}$.

Theorem \ref{thm:InvWeylDn} says that, over a field $k_0$ whose characteristic does not divide the order of $\Weyl(D_n)$, $\Inv_{k_0}(\Weyl(D_n), W)$ is generated as a $W(k_0)$-algebra by trace form invariants, namely $a_{L/k}'$, $a_{K/k}'$, and $r$, and their exterior powers. 

\begin{theorem}\label{thm:InvWeylDn}
Suppose that $\textup{char}(k_0)$ does not divide $|\Weyl(D_n)|$, that is, $\textup{char}(k_0)$ is zero or greater than $n$. If $n$ is odd, then $\Inv_{k_0}(\Weyl(D_n), W)$ is generated as a $W(k_0)$-algebra by $a'_{L/k}$, $a'_{K/k}$, and their exterior powers. If $n$ is even, then $\Inv_{k_0}(\Weyl(D_n), W)$ is generated as a $W(k_0)$-algebra by $a'_{L/k}$, $a'_{K/k}$, their exterior powers, and $r$.
\end{theorem}

The proof of Theorem \ref{thm:InvWeylDn} is closely related to the proof of Theorem \ref{thm:InvWeylBn}.

\begin{proof}[Proof of Theorem \ref{thm:InvWeylDn}.]

In \cite[Corollary 7.5, p\ 799]{hirsch2020decomposability}, Hirsch describes generators\footnote{Note two points of possible confusion here. First, the invariants that we refer to as $u_d \circ \varphi_*$ and $v_d' \circ \varphi_*$ are referred to in \cite[Section 7, p\ 793-799]{hirsch2020decomposability} as $u_d$ and $v_d'$ respectively. We write $\varphi_*$ here to distinguish between invariants of $\Weyl(D_n)$ and invariants of $\Weyl(B_n)$. Second, Hirsch's definition of the invariant $v_d'$ of $\Weyl(D_n)$ is not completely explicit. One can check that this invariant is indeed the one we call $v_d' \circ \varphi_*$, by comparing the invariant $\res_{\Weyl(B_n)}^{P_L}(v_d')$ \cite[Lemma 5.3, p\ 786]{hirsch2020decomposability} for $L = \floor{n/2}$ with the invariant $\res_{\Weyl(D_n)}^P(v_d')$ \cite[p\ 795]{hirsch2020decomposability} and applying the splitting principle \cite[Proposition 2.3, p\ 770]{hirsch2020decomposability}.} for the $H(k_0)$-module $\Inv_{k_0}(\Weyl(D_n), H)$, for $n \geq 2$. These module generators are contained in the $H(k_0)$-subalgebra of $\Inv_{k_0}(\Weyl(D_n), H)$ generated by
\[
	\{u_d \circ \varphi_*\}_{0 \leq d \leq n} \cup \{v_d' \circ \varphi_*\}_{0 \leq d \leq n} \cup R,
\]
where $u_d$ and $v_d'$ are the Witt invariants of $\Weyl(B_n)$ described in Section \ref{subsec:winv_Bn}, and $R$ is defined by $R = \emptyset$ if $n$ is odd and $R = e_{n/2} \circ r$ if $n$ is even. The invariant $e_{n/2} \circ r$ is well-defined by \cite[(7.2), p\ 796]{hirsch2020decomposability}, which shows that $r$ maps $H^1(k, \Weyl(D_n))$ into $I^{n/2}(k)$.

By the proof of Theorem \ref{thm:InvWeylBn}, the invariants $u_d$ and $v_d'$ lift to Witt invariants in the $W(k_0)$-subalgebra of $\Inv_{k_0}(\Weyl(B_n), W)$ generated by the exterior powers of $a_{L/k}$ and $a_{K/k}$. So $u_d \circ \varphi_*$ and $v_d' \circ \varphi_*$ lift to elements of $\Inv_{k_0}(\Weyl(D_n), W)$ in the $W(k_0)$-algebra generated by the exterior powers of $a_{L/k} \circ \varphi_*$ and $a_{K/k} \circ \varphi_*$. By Lemma \ref{lem:Weyl_Char_Implies_BddInvar} and Proposition \ref{prop:InvGenMod}, $\Inv_{k_0}(\Weyl(D_n), W)$ is generated as a $W(k_0)$-algebra by the exterior powers of $a_{L/k} \circ \varphi_*$ and $a_{K/k} \circ \varphi_*$, together with $r$ when $n$ is even.
\end{proof}

\begin{remark}\label{rem:OtherExceptionalTypes}
It would be interesting to extend the techniques in Section \ref{subsec:winv_Bn} to describe the Witt invariants of the Weyl groups of the remaining types, namely $F_4$, $E_6$, $E_7$, and $E_8$. For such a Weyl group $G$, is $\Inv_{k_0}(G, W)$ generated as a $W(k_0)$-algebra by trace forms and their exterior powers? 

Following the proof for type $B_n$ (Theorem \ref{thm:InvWeylBn}), one would hope to find generators for $\Inv_{k_0}(G, H)$ of the form
\[
	H^1(-, G) \to H^1(-, S_n) \xrightarrow{(\iota_{n})_*} H^1(-, O_n) \xrightarrow{\swmod_d} H^d(-, \Ztwo),
\]
for various $n$ and $d$ and maps $H^1(-, G) \to H^1(-, S_n)$. Unlike the $B_n$ case, however, the existence of such generators is not immediate from \cite{hirsch2020decomposability}. For example, for $G$ the Weyl group of type $F_4$, Hirsch describes certain generators of $\Inv_{k_0}(G, H)$ by taking the inclusion of $G$ into $O_4$ as an orthogonal reflection group, then composing with Stiefel-Whitney classes \cite[p\ 792-793]{hirsch2020decomposability}. The embedding $G \to O_4$ does not factor through $S_4$, so we cannot directly apply the reasoning of Section \ref{subsec:winv_Bn}. 
\end{remark}

\bibliography{bibliography}

\providecommand{\bysame}{\leavevmode\hbox to3em{\hrulefill}\thinspace}
\providecommand{\MR}{\relax\ifhmode\unskip\space\fi MR }
\providecommand{\MRhref}[2]{%
  \href{http://www.ams.org/mathscinet-getitem?mr=#1}{#2}
}
\providecommand{\href}[2]{#2}
\begin{thebibliography}{KMRT98}

\bibitem[AP71]{ArasonPfister1971}
J{\'o}n~Kristinn Arason and Albrecht Pfister, \emph{Beweis des {K}rullschen
  {D}urchschnittsatzes f{\"u}r den {W}ittring}, Inventiones mathematicae
  \textbf{12} (1971), no.~2, 173--176.

\bibitem[Bou07]{Bourbaki2007}
Nicolas Bourbaki, \emph{Alg{\`e}bre: Chapitre 9}, Springer, Berlin, Heidelberg,
  2007.

\bibitem[CP84]{conner1984survey}
P.E. Conner and R.V. Perlis, \emph{A survey of trace forms of algebraic number
  fields}, Series In Pure Mathematics, World Scientific Publishing Company,
  1984.

\bibitem[GH21]{GILLE2021}
Stefan Gille and Christian Hirsch, \emph{On the splitting principle for
  cohomological invariants of reflection groups}, Transformation Groups (2021).

\bibitem[GMS03]{GaribaldiMerkurjevSerre03}
Skip Garibaldi, Alexander Merkurjev, and Jean-Pierre Serre, \emph{Cohomological
  invariants in {Galois} cohomology}, University Lecture Series, vol.~28,
  American Mathematical Society, Providence, Rhode Island, 2003.

\bibitem[GS17]{gillesamuely2017central}
P.~Gille and T.~Szamuely, \emph{Central simple algebras and {G}alois
  cohomology}, Cambridge Studies in Advanced Mathematics, Cambridge University
  Press, 2017.

\bibitem[Hir20]{hirsch2020decomposability}
Christian Hirsch, \emph{On the decomposability of mod 2 cohomological
  invariants of {Weyl} groups}, Comment. Math. Helv. \textbf{95} (2020),
  765--–809.

\bibitem[Hum72]{Humphreys1972}
James~E. Humphreys, \emph{Introduction to {L}ie algebras and representation
  theory}, Springer New York, New York, NY, 1972.

\bibitem[KMRT98]{KMRT1998}
Max-Albert Knus, Alexander Merkurjev, Markus Rost, and Jean-Pierre Tignol,
  \emph{The book of involutions}, Colloquium Publications, vol.~44, American
  Mathematical Society, Providence, Rhode Island, 1998.

\bibitem[Mor05]{Morel2005Milnors}
Fabien Morel, \emph{Milnor’s conjecture on quadratic forms and mod 2 motivic
  complexes}, Rendiconti del Seminario Matematico della Università di Padova
  \textbf{114} (2005).

\bibitem[Neu99]{Neukirch1999}
J{\"u}rgen Neukirch, \emph{Algebraic number theory}, Springer Berlin
  Heidelberg, Berlin, Heidelberg, 1999.

\bibitem[OVV07]{OrlovVishikVoevodsky2007}
D.~Orlov, A.~Vishik, and V.~Voevodsky, \emph{An exact sequence for
  {{\(K^M_*/2\)}} with applications to quadratic forms}, Ann. Math. (2)
  \textbf{165} (2007), no.~1, 1--13 (English).

\bibitem[Ros96]{Rost1996}
Markus Rost, \emph{Chow groups with coefficients.}, Documenta Mathematica
  \textbf{1} (1996), 319--393 (eng).

\bibitem[Ser97]{Serre1997}
Jean-Pierre Serre, \emph{{Galois} cohomology}, Springer, Berlin, Heidelberg,
  1997.

\bibitem[Ser03]{Serre03_GMS}
\bysame, \emph{Cohomological invariants, {W}itt invariants, and trace forms},
  Cohomological invariants in {G}alois cohomology, Univ. Lecture Ser., vol.~28,
  Amer. Math. Soc., Providence, RI, 2003, Notes by Skip Garibaldi, pp.~1--100.

\bibitem[Ser18]{serre2018cohomological}
\bysame, \emph{Cohomological invariants mod 2 of {Weyl} groups}, Oberwolfach
  Rep., vol.~21, 2018, pp.~1284--1286.

\bibitem[Spr98]{springer1998linear}
T.A. Springer, \emph{Linear algebraic groups}, Progress in mathematics,
  Springer, 1998.

\bibitem[Voe03]{Voevodsky2003}
Vladimir Voevodsky, \emph{Motivic cohomology with $\mathbb{Z}/2$-coefficients},
  Publications Math{\'e}matiques de l'Institut des Hautes {\'E}tudes
  Scientifiques \textbf{98} (2003), no.~1, 59--104.

\end{thebibliography}
\bibliographystyle{amsalpha}

\end{document}